\documentclass{article}
\usepackage[utf8]{inputenc}
\usepackage[letterpaper,margin=1in]{geometry}
\usepackage{amsmath, amssymb, amsthm, mathtools}
\usepackage{enumitem}

\newcommand{\ds}{\displaystyle}
\newtheorem{theorem}{Theorem}[section]
\newtheorem{lemma}[theorem]{Lemma}
\newtheorem{corollary}[theorem]{Corollary}
\theoremstyle{definition}\newtheorem{remark}[theorem]{Remark}

\DeclareMathOperator{\Child}{Chi}
\DeclareMathOperator{\parent}{par}
\DeclareMathOperator{\Gen}{Gen}
\newcommand{\Z}{\mathbb{Z}}
\newcommand{\R}{\mathbb{R}}
\newcommand{\N}{\mathbb{N}}
\newcommand{\C}{\mathbb{C}}
\newcommand{\K}{\mathbb{K}}
\begin{document}

\title{Positive expansivity of weighted forward shifts on rooted directed trees}

\author{
  Sonja Ebey, Keri Hoatson, and David Walmsley\thanks{Corresponding author: \texttt{walmsl1@stolaf.edu}} \\[1ex]
  \small Department of Mathematics, Statistics, and Computer Science \\
  \small St.\ Olaf College, Northfield, MN 55057, USA
}

\date{\today}
\maketitle

\begin{abstract}
    We characterize, in terms of the weights, when a weighted forward shift on the space $\ell^p(V)$, $1\leq p<\infty$, or $c_0(V)$ of a rooted directed tree $V$ is positively expansive, average positively expansive, or uniformly positively expansive. For symmetric weights on symmetric trees, the conditions reduce to the corresponding conditions for weighted forward shifts on $\N_0$ with suitably modified weights.
\end{abstract}

\noindent\textit{2020 Mathematics Subject Classification:} 47A16, 47B37, 37B05.

\noindent\textit{Keywords:} Positive expansivity, average positive expansivity, uniform positive expansivity, weighted forward shift, directed tree.

\section{Introduction}\label{Introduction}

Expansivity is one of the basic notions in the theory of dynamical systems. A homeomorphism $f$ of a metric space $(M,d)$ is \emph{expansive} if there is a constant $c>0$ such that for any two distinct points $x,y\in M$ there exists $n\in\Z$ with $d(f^n(x),f^n(y))>c$. In linear dynamics, an invertible bounded linear operator $T$ on a Banach space $X$ is called \emph{expansive} if for every $x$ in the unit sphere $S_X$ of $X$ there exists $n\in\Z$ with $\|T^nx\|\geq2$, and \emph{uniformly expansive} if there exists $n\in\N$ such that $\|T^nx\|\geq2$ or $\|T^{-n}x\|\geq2$ for every $x\in S_X$. Eisenberg \cite{Eisenberg} showed that an invertible operator on $\C^d$ is expansive if and only if it has no eigenvalue of modulus one, and Eisenberg and Hedlund \cite{EisenbergHedlund} and Hedlund \cite{Hedlund} showed that an invertible operator on a Banach space is uniformly expansive if and only if its approximate point spectrum does not meet the unit circle. These notions were revisited by Bernardes, Cirilo, Darji, Messaoudi and Pujals \cite{(9)Bernardes}, who studied their relations with hyperbolicity and shadowing and characterized them for weighted shifts on the classical sequence spaces.

The definitions above involve negative powers of $T$ and therefore require $T$ to be invertible. Their \emph{positive} versions involve only forward iterates and make sense for every bounded linear operator $T$ on $X$. Following \cite{(9)Bernardes} and \cite{(1)Alves}, $T$ is called
\begin{itemize}
    \item \emph{positively expansive} if for every $x\in S_X$ there exists $n\in\N$ with $\|T^nx\|\geq 2$;
    \item \emph{uniformly positively expansive} if there exists $n\in\N$ such that $\|T^nx\|\geq 2$ for every $x\in S_X$;
    \item \emph{average positively expansive} if $\ds\sup_{n\in\N}\frac1n\sum_{j=0}^{n-1}\|T^jx\|=\infty$ for every nonzero $x\in X$.
\end{itemize}
By \cite[Proposition 19]{(9)Bernardes}, $T$ is positively expansive if and only if $\sup_{n\in\N}\|T^nx\|=\infty$ for every nonzero $x\in X$, and uniformly positively expansive if and only if $\|T^nx\|\to\infty$ uniformly on $S_X$. Average positive expansivity was introduced by Alves, Bernardes and Messaoudi \cite{(1)Alves} as the positive counterpart of average expansivity, a notion they defined in connection with the average shadowing property; among other results, they showed that an invertible operator is hyperbolic if and only if it is average expansive and has the average shadowing property \cite[Theorem 14]{(1)Alves}. The three positive notions are related by
\[
    \text{uniform positive expansivity}\ \Longrightarrow\ \text{average positive expansivity}\ \Longrightarrow\ \text{positive expansivity},
\]
and none of the converse implications holds in general \cite[Remark 18]{(1)Alves}. A unilateral weighted backward shift is never positively expansive, while the unilateral weighted forward shifts with each of the three properties are characterized by conditions on the products of the weights in \cite[Remark 32]{(9)Bernardes} and \cite[Remark 18]{(1)Alves}.

Weighted shifts are a central source of examples in linear dynamics; see \cite{BayartMatheron,GrossePeris}. A natural way to enlarge this class is to replace the index sets $\N$ and $\Z$ by more general directed trees. Weighted shifts on directed trees were introduced by Jab\l o\'nski, Jung and Stochel \cite{Jablonski}, who studied their operator-theoretic properties on $\ell^2$-spaces, and the study of their dynamics was initiated by Mart\'inez-Avenda\~no \cite{Martinez}. Grosse-Erdmann and Papathanasiou \cite{GrossePapahypercyclic} characterized the bounded weighted shifts on the spaces $\ell^p(V)$, $1\leq p<\infty$, and $c_0(V)$ of a directed tree $V$, as well as the hypercyclic, weakly mixing and mixing weighted backward shifts on these spaces. Subsequent work characterized, among others, the chaotic \cite{GrossePapachaotic}, the $\mathcal{F}$-transitive and $\mathcal{F}$-recurrent \cite{AbakumovAbbar}, and the chain recurrent \cite{MortensenWalmsley} weighted backward shifts on trees; see also \cite{LopezPapachainrecontrees}.

Weighted forward shifts on rooted trees behave very differently. Since every function in the range of a forward shift vanishes at the root, such a shift is never hypercyclic \cite{Martinez,GrossePapahypercyclic}, nor even chain recurrent (see \cite[Remark 21]{(1)Alves} for $V=\N_0$, and \cite[Proposition 2.6]{MortensenWalmsley}). Instead, the orbits of a forward shift spread out through the tree, and it is natural to ask when they grow. In this note we answer this question for the three positive expansivity notions above: we characterize, in terms of the weights, when a weighted forward shift on $\ell^p(V)$, $1\leq p<\infty$, or on $c_0(V)$ of a rooted directed tree $V$ is positively expansive, average positively expansive, or uniformly positively expansive (Theorems \ref{main theorem} and \ref{main theorem c0}). Our approach rests on an explicit formula for the norms of the iterates of a forward shift (Lemma \ref{norm lemma}). When both the tree and the weight are symmetric, the conditions reduce to the known conditions for weighted forward shifts on $\N_0$ with suitably modified weights (Corollary \ref{symmetric corollary}).

\section{The characterization}\label{Characterization section}

We use the terminology of \cite{Jablonski,GrossePapahypercyclic}. A \emph{directed tree} $(V,E)$ is a connected directed graph, with a countable set $V$ of vertices and a set $E\subset(V\times V)\setminus\{(v,v):v\in V\}$ of edges, that has no cycles and in which each vertex $v$ has at most one \emph{parent}, that is, a vertex $\parent(v)$ with $(\parent(v),v)\in E$. At most one vertex has no parent; if it exists, it is called the \emph{root}, denoted by $\mathsf{root}$, and the tree is called \emph{rooted}. Throughout, $(V,E)$ is a rooted directed tree. The set of \emph{children} of $v\in V$ is $\Child(v)=\{u\in V:\parent(u)=v\}$, and a vertex without children is a \emph{leaf}. We let $\Child^0(v)=\{v\}$ and, for $n\geq1$, we let $\Child^n(v)$ be the set of children of the vertices in $\Child^{n-1}(v)$; thus $u\in\Child^n(v)$ if and only if $\parent^n(u)=v$. The \emph{generations} $\Gen_n=\Child^n(\mathsf{root})$, $n\in\N_0$, partition $V$.

Let $\K=\R$ or $\C$. For $1\leq p<\infty$, $\ell^p(V)$ is the space of all $f\in\K^V$ with $\|f\|=\big(\sum_{v\in V}|f(v)|^p\big)^{1/p}<\infty$, and $c_0(V)$ is the space of all $f\in\K^V$ for which $\{v\in V:|f(v)|\geq\varepsilon\}$ is finite for every $\varepsilon>0$, endowed with the norm $\|f\|=\sup_{v\in V}|f(v)|$. Throughout, $X$ denotes one of these spaces, and $e_v$, $v\in V$, denote the canonical unit vectors.

A \emph{weight} is a family $\lambda=(\lambda_v)_{v\in V}$ of nonzero scalars. The \emph{weighted forward shift} $S_\lambda$ is defined by
\[
    [S_\lambda(f)](v)=\lambda_vf(\parent(v))\quad\text{for } v\neq\mathsf{root},\qquad [S_\lambda(f)](\mathsf{root})=0.
\]
By \cite[Proposition 2.2]{GrossePapahypercyclic}, $S_\lambda$ is a bounded operator on $\ell^p(V)$ if and only if $\sup_{v\in V}\sum_{u\in\Child(v)}|\lambda_u|^p<\infty$, and on $c_0(V)$ if and only if $\sup_{v\neq\mathsf{root}}|\lambda_v|<\infty$ and $(\lambda_u)_{u\in\Child(v)}\in c_0(\Child(v))$ for every $v\in V$. As in \cite{GrossePapahypercyclic}, for $v\in V$ and $u\in\Child^n(v)$ with $n\geq1$ we write
\[
    \lambda(v\to u)\coloneqq\lambda_{\parent^{n-1}(u)}\cdots\lambda_{\parent(u)}\lambda_u
\]
for the product of the weights $\lambda_w$ over the vertices $w\neq v$ on the path from $v$ to $u$, and we set $\lambda(v\to v)\coloneqq1$. Then $\lambda(v\to u)=\lambda(v\to w)\lambda(w\to u)$ for every vertex $w$ on the path from $v$ to $u$. A sum or a supremum over the empty set is taken to be $0$.

\begin{lemma}\label{norm lemma}
    Let $X=\ell^p(V)$, $1\leq p<\infty$, or $X=c_0(V)$, let $\lambda$ be a weight for which $S_\lambda$ is bounded on $X$, and let $n\in\N_0$. For every $f\in X$ and $u\in V$,
    \begin{equation}\label{iterate formula}
        [S_\lambda^n(f)](u)=
        \begin{cases}
            \lambda(\parent^n(u)\to u)\,f(\parent^n(u)) & \text{if } \parent^n(u) \text{ exists},\\
            0 & \text{otherwise},
        \end{cases}
    \end{equation}
    and
    \begin{align}
        \|S_\lambda^n(f)\|^p&=\sum_{v\in V}|f(v)|^p\sum_{u\in\Child^n(v)}|\lambda(v\to u)|^p &&\text{if } X=\ell^p(V),\label{norm formula lp}\\
        \|S_\lambda^n(f)\|&=\sup_{v\in V}|f(v)|\sup_{u\in\Child^n(v)}|\lambda(v\to u)| &&\text{if } X=c_0(V).\label{norm formula c0}
    \end{align}
    In particular, $\|S_\lambda^n(f)\|\geq|f(v)|\,\|S_\lambda^n(e_v)\|$ for all $f\in X$ and $v\in V$, and $\inf_{f\in S_X}\|S_\lambda^n(f)\|=\inf_{v\in V}\|S_\lambda^n(e_v)\|$.
\end{lemma}

\begin{proof}
    Formula \eqref{iterate formula} holds for $n=0$. If it holds for some $n$, then $[S_\lambda^{n+1}(f)](\mathsf{root})=0$, while for $u\neq\mathsf{root}$,
    \begin{align*}
        [S_\lambda^{n+1}(f)](u)&=\lambda_u\,[S_\lambda^n(f)](\parent(u))=\lambda_u\,\lambda(\parent^{n+1}(u)\to\parent(u))\,f(\parent^{n+1}(u))\\
        &=\lambda(\parent^{n+1}(u)\to u)\,f(\parent^{n+1}(u))
    \end{align*}
    if $\parent^{n+1}(u)$ exists, and $[S_\lambda^{n+1}(f)](u)=0$ otherwise. This proves \eqref{iterate formula}. Since each vertex has at most one $n$th ancestor, the sets $\Child^n(v)$, $v\in V$, are pairwise disjoint, and their union is the set of vertices $u$ for which $\parent^n(u)$ exists. Grouping these vertices according to $v=\parent^n(u)$ in \eqref{iterate formula} gives \eqref{norm formula lp} and, with suprema in place of sums, \eqref{norm formula c0}. Taking $f=e_v$ shows that the inner sum in \eqref{norm formula lp} equals $\|S_\lambda^n(e_v)\|^p$ and that the inner supremum in \eqref{norm formula c0} equals $\|S_\lambda^n(e_v)\|$. Hence $\|S_\lambda^n(f)\|^p=\sum_{v\in V}|f(v)|^p\,\|S_\lambda^n(e_v)\|^p$, respectively $\|S_\lambda^n(f)\|=\sup_{v\in V}|f(v)|\,\|S_\lambda^n(e_v)\|$, and the remaining assertions follow; for the last one, note that $\|S_\lambda^n(f)\|\geq\inf_{w\in V}\|S_\lambda^n(e_w)\|\,\|f\|$ for every $f\in X$.
\end{proof}

By Lemma \ref{norm lemma}, the orbit of a function $f\in X$ dominates, up to the factor $|f(v)|$, the orbit of $e_v$ for every vertex $v$ in the support of $f$. This reduces the three expansivity properties of $S_\lambda$ to conditions on the canonical vectors, and hence on the weights.

\begin{theorem}\label{main theorem}
    Let $(V,E)$ be a rooted directed tree, let $1\leq p<\infty$, and let $\lambda$ be a weight for which $S_\lambda$ is bounded on $\ell^p(V)$. Then:
    \begin{enumerate}[label=\upshape(\alph*)]
        \item $S_\lambda$ is positively expansive on $\ell^p(V)$ if and only if, for every $v\in V$,
        \[
            \sup_{n\in\N}\,\sum_{u\in\Child^n(v)}|\lambda(v\to u)|^p=\infty;
        \]
        \item $S_\lambda$ is average positively expansive on $\ell^p(V)$ if and only if, for every $v\in V$,
        \[
            \sup_{n\in\N}\,\frac1n\sum_{j=0}^{n-1}\Big(\sum_{u\in\Child^j(v)}|\lambda(v\to u)|^p\Big)^{1/p}=\infty;
        \]
        \item the following assertions are equivalent:
        \begin{enumerate}[label=\upshape(\roman*)]
            \item $S_\lambda$ is uniformly positively expansive on $\ell^p(V)$;
            \item there exists $N\in\N$ such that $\ds\inf_{v\in V}\,\sum_{u\in\Child^N(v)}|\lambda(v\to u)|^p>1$;
            \item $\ds\lim_{n\to\infty}\,\inf_{v\in V}\,\sum_{u\in\Child^n(v)}|\lambda(v\to u)|^p=\infty$.
        \end{enumerate}
    \end{enumerate}
\end{theorem}

\begin{theorem}\label{main theorem c0}
    Let $(V,E)$ be a rooted directed tree, and let $\lambda$ be a weight for which $S_\lambda$ is bounded on $c_0(V)$. Then:
    \begin{enumerate}[label=\upshape(\alph*)]
        \item $S_\lambda$ is positively expansive on $c_0(V)$ if and only if, for every $v\in V$,
        \[
            \sup_{n\in\N}\,\sup_{u\in\Child^n(v)}|\lambda(v\to u)|=\infty;
        \]
        \item $S_\lambda$ is average positively expansive on $c_0(V)$ if and only if, for every $v\in V$,
        \[
            \sup_{n\in\N}\,\frac1n\sum_{j=0}^{n-1}\,\sup_{u\in\Child^j(v)}|\lambda(v\to u)|=\infty;
        \]
        \item the following assertions are equivalent:
        \begin{enumerate}[label=\upshape(\roman*)]
            \item $S_\lambda$ is uniformly positively expansive on $c_0(V)$;
            \item there exists $N\in\N$ such that $\ds\inf_{v\in V}\,\sup_{u\in\Child^N(v)}|\lambda(v\to u)|>1$;
            \item $\ds\lim_{n\to\infty}\,\inf_{v\in V}\,\sup_{u\in\Child^n(v)}|\lambda(v\to u)|=\infty$.
        \end{enumerate}
    \end{enumerate}
\end{theorem}

\begin{proof}[Proof of Theorems \ref{main theorem} and \ref{main theorem c0}]
    Let $X=\ell^p(V)$ in the case of Theorem \ref{main theorem} and $X=c_0(V)$ in the case of Theorem \ref{main theorem c0}. By Lemma \ref{norm lemma}, for all $v\in V$ and $n\in\N_0$,
    \[
        \|S_\lambda^n(e_v)\|=\Big(\sum_{u\in\Child^n(v)}|\lambda(v\to u)|^p\Big)^{1/p},\quad\text{respectively}\quad \|S_\lambda^n(e_v)\|=\sup_{u\in\Child^n(v)}|\lambda(v\to u)|.
    \]
    Since $t\mapsto t^p$ is increasing on $[0,\infty)$, the conditions in (a), (b), (c)(ii) and (c)(iii) of both theorems therefore say, respectively, that $\sup_{n\in\N}\|S_\lambda^n(e_v)\|=\infty$ for every $v\in V$, that $\sup_{n\in\N}\frac1n\sum_{j=0}^{n-1}\|S_\lambda^j(e_v)\|=\infty$ for every $v\in V$, that $\inf_{v\in V}\|S_\lambda^N(e_v)\|>1$ for some $N\in\N$, and that $\lim_{n\to\infty}\inf_{v\in V}\|S_\lambda^n(e_v)\|=\infty$.

    (a) If $S_\lambda$ is positively expansive, then $\sup_{n}\|S_\lambda^n(e_v)\|=\infty$ for every $v\in V$ by \cite[Proposition 19]{(9)Bernardes}. Conversely, assume that this holds, and let $x\in S_X$. Choose $v\in V$ with $x(v)\neq0$. By Lemma \ref{norm lemma}, $\|S_\lambda^n(x)\|\geq|x(v)|\,\|S_\lambda^n(e_v)\|$ for every $n$, so $\|S_\lambda^n(x)\|\geq2$ for some $n\in\N$.

    (b) If $S_\lambda$ is average positively expansive, the condition follows by taking $x=e_v$ in the definition. Conversely, let $x\in X$ be nonzero, and choose $v\in V$ with $x(v)\neq0$. By Lemma \ref{norm lemma},
    \[
        \frac1n\sum_{j=0}^{n-1}\|S_\lambda^j(x)\|\geq|x(v)|\,\frac1n\sum_{j=0}^{n-1}\|S_\lambda^j(e_v)\|\quad\text{for every } n\in\N,
    \]
    and the right-hand side is unbounded by assumption.

    (c) Let $m_n\coloneqq\inf_{v\in V}\|S_\lambda^n(e_v)\|$, so that $m_n=\inf_{x\in S_X}\|S_\lambda^n(x)\|$ by Lemma \ref{norm lemma}. Thus (i) holds if and only if $m_N\geq2$ for some $N\in\N$. Moreover, $m_{j+k}\geq m_jm_k$ for all $j,k\in\N_0$, since $\|S_\lambda^{j+k}(x)\|\geq m_j\|S_\lambda^k(x)\|\geq m_jm_k\|x\|$ for every $x\in X$. Clearly (i) implies (ii), and (iii) implies (i). Suppose that (ii) holds. Since $\|S_\lambda^N(x)\|\leq\|S_\lambda\|^{N-1}\|S_\lambda(x)\|$ for every $x\in X$, we have $0<m_N\leq\|S_\lambda\|^{N-1}m_1$, so $m_1>0$. Writing $n=kN+r$ with $k\in\N_0$ and $0\leq r<N$, we obtain
    \[
        m_n\geq m_N^k\,m_1^r\geq m_N^k\min\{1,m_1\}^{N-1}\to\infty\quad\text{as } n\to\infty,
    \]
    which is (iii).
\end{proof}

\begin{remark}\label{main remarks}
    \begin{enumerate}[label=\upshape(\alph*)]
        \item The conditions in parts (a) and (b) of Theorems \ref{main theorem} and \ref{main theorem c0} are needed at every vertex; it does not suffice to check them at the root. For example (compare the discussion following \cite[Example 4.6]{GrossePapahypercyclic}), let the root have two children $a_1$ and $b_1$, and let every other vertex have exactly one child, so that $V\setminus\{\mathsf{root}\}$ consists of two rays $a_1\to a_2\to\cdots$ and $b_1\to b_2\to\cdots$. If $\lambda_{a_i}=2$ and $\lambda_{b_i}=1$ for all $i$, then $\lambda(\mathsf{root}\to a_n)=2^n$, so these conditions hold at the root, while $\Child^n(b_1)=\{b_{n+1}\}$ and $\lambda(b_1\to b_{n+1})=1$ for every $n$, so they fail at $b_1$. Thus $S_\lambda$ is not positively expansive on any of the spaces $\ell^p(V)$ or $c_0(V)$.
        \item If $V$ has a leaf $v$, then $\Child^n(v)=\emptyset$ for every $n\geq1$, so the sums and suprema in Theorems \ref{main theorem} and \ref{main theorem c0} vanish at $v$, and $S_\lambda$ has none of the three properties; this is also clear from $S_\lambda(e_v)=0$.
    \end{enumerate}
\end{remark}

We now turn to symmetric weighted forward shifts on symmetric trees, the setting closest to the classical shifts. A rooted directed tree is \emph{symmetric} if any two vertices in the same generation have the same number of children; for $v\in\Gen_k$ we then write $\gamma_k=|\Child(v)|$. A weight $\lambda$ is \emph{symmetric} if $\lambda_u=\lambda_v$ whenever $u$ and $v$ lie in the same generation; for $v\in\Gen_k$, $k\geq1$, we then write $\lambda_k=\lambda_v$ (the weight of the root plays no role). As for the backward shifts in \cite[Corollary 4.7]{GrossePapahypercyclic}, it then suffices to check the conditions in parts (a) and (b) of Theorems \ref{main theorem} and \ref{main theorem c0} at the root.

\begin{corollary}\label{symmetric corollary}
    Let $(V,E)$ be a symmetric rooted directed tree without leaves, let $X=\ell^p(V)$, $1\leq p<\infty$, or $X=c_0(V)$, and let $\lambda=(\lambda_k)_{k\geq1}$ be a symmetric weight for which $S_\lambda$ is bounded on $X$. For $k\geq1$, let
    \[
        w_k\coloneqq\gamma_{k-1}^{1/p}\,|\lambda_k| \ \text{ if } X=\ell^p(V),\qquad w_k\coloneqq|\lambda_k| \ \text{ if } X=c_0(V).
    \]
    Then:
    \begin{enumerate}[label=\upshape(\alph*)]
        \item $S_\lambda$ is positively expansive if and only if $\ds\sup_{n\in\N}\,w_1\cdots w_n=\infty$;
        \item $S_\lambda$ is average positively expansive if and only if $\ds\sup_{n\in\N}\frac1n\sum_{j=1}^{n}w_1\cdots w_j=\infty$;
        \item $S_\lambda$ is uniformly positively expansive if and only if $\ds\lim_{n\to\infty}\,\inf_{k\in\N_0}w_{k+1}\cdots w_{k+n}=\infty$.
    \end{enumerate}
\end{corollary}

\begin{proof}
    Note that each $\gamma_k$ is finite, since $S_\lambda$ is bounded. Let $v\in\Gen_k$, $k\in\N_0$, and $n\in\N_0$. Every $u\in\Child^n(v)$ lies in $\Gen_{k+n}$, so that $|\lambda(v\to u)|=|\lambda_{k+1}\cdots\lambda_{k+n}|$; moreover, $|\Child^n(v)|=\gamma_k\cdots\gamma_{k+n-1}\geq1$. Hence
    \begin{align*}
        \Big(\sum_{u\in\Child^n(v)}|\lambda(v\to u)|^p\Big)^{1/p}&=w_{k+1}\cdots w_{k+n} &&\text{if } X=\ell^p(V),\\
        \sup_{u\in\Child^n(v)}|\lambda(v\to u)|&=w_{k+1}\cdots w_{k+n} &&\text{if } X=c_0(V),
    \end{align*}
    and in both cases $w_{k+1}\cdots w_{k+n}=P_{k+n}/P_k$, where $P_m\coloneqq w_1\cdots w_m$ and $P_0\coloneqq1$. Part (c) now follows from part (c) of Theorems \ref{main theorem} and \ref{main theorem c0}, after taking $p$th roots in the case of $\ell^p(V)$. Since $\sup_nP_{k+n}/P_k=\infty$ if and only if $\sup_nP_n=\infty$, part (a) follows in the same way from part (a) of these theorems. For (b), let $c_n\coloneqq\frac1n\sum_{j=0}^{n-1}P_j$. Since $c_n\leq1+\frac1n\sum_{j=1}^{n}P_j$ and $c_{n+1}\geq\frac12\cdot\frac1n\sum_{j=1}^nP_j$, the condition in (b) is equivalent to $\sup_nc_n=\infty$, which is the condition in part (b) of the theorems at the root. If $\sup_nc_n=\infty$, then for every $v\in\Gen_k$,
    \[
        \frac1n\sum_{j=0}^{n-1}\frac{P_{k+j}}{P_k}=\frac{1}{P_k}\Big(\frac{n+k}{n}\,c_{n+k}-\frac1n\sum_{i=0}^{k-1}P_i\Big)\geq\frac{1}{P_k}\Big(c_{n+k}-\sum_{i=0}^{k-1}P_i\Big),
    \]
    and the right-hand side is unbounded in $n$, because each $c_m$ is finite and therefore $\sup_{m>k}c_m=\infty$. Thus the condition in part (b) of the theorems holds at every vertex, and (b) follows.
\end{proof}

As the proof shows, the symmetry of the tree is not needed when $X=c_0(V)$. For $X=\ell^1(V)$, the products in Corollary \ref{symmetric corollary} are $w_{k+1}\cdots w_{k+n}=|\gamma_k\lambda_{k+1}\gamma_{k+1}\lambda_{k+2}\cdots\gamma_{k+n-1}\lambda_{k+n}|$. For the tree $V=\N_0$, in which every vertex has exactly one child, we have $\gamma_k=1$ for all $k$, and Corollary \ref{symmetric corollary} recovers the characterizations of the unilateral weighted forward shifts with the three properties in \cite[Remark 32]{(9)Bernardes} and \cite[Remark 18]{(1)Alves}.

\bibliography{expansivity}

\begin{thebibliography}{10}

\bibitem{AbakumovAbbar}
Evgeny Abakumov and Arafat Abbar.
\newblock On several dynamical properties of shifts acting on directed trees.
\newblock {\em Canad. J. Math.}, 78(2):375--404, 2026.

\bibitem{(1)Alves}
Fabricio~F. Alves, Nilson~C. Bernardes, Jr., and Ali Messaoudi.
\newblock Chain recurrence and average shadowing in dynamics.
\newblock {\em Monatsh. Math.}, 196(4):665--697, 2021.

\bibitem{BayartMatheron}
Fr\'{e}d\'{e}ric Bayart and \'{E}tienne Matheron.
\newblock {\em Dynamics of linear operators}, volume 179 of {\em Cambridge Tracts in Mathematics}.
\newblock Cambridge University Press, Cambridge, 2009.

\bibitem{(9)Bernardes}
Nilson~C. Bernardes, Jr., Patricia~R. Cirilo, Udayan~B. Darji, Ali Messaoudi, and Enrique~R. Pujals.
\newblock Expansivity and shadowing in linear dynamics.
\newblock {\em J. Math. Anal. Appl.}, 461(1):796--816, 2018.

\bibitem{Eisenberg}
M.~Eisenberg.
\newblock Expansive automorphisms of finite dimensional vector spaces.
\newblock {\em Fund. Math.}, 59:307--312, 1966.

\bibitem{EisenbergHedlund}
M.~Eisenberg and J.~H. Hedlund.
\newblock Expansive automorphisms of {B}anach spaces.
\newblock {\em Pacific J. Math.}, 34(3):647--656, 1970.

\bibitem{GrossePapahypercyclic}
Karl-G. Grosse-Erdmann and Dimitris Papathanasiou.
\newblock Dynamics of weighted shifts on directed trees.
\newblock {\em Indiana Univ. Math. J.}, 72(1):263--299, 2023.

\bibitem{GrossePapachaotic}
Karl-G. Grosse-Erdmann and Dimitris Papathanasiou.
\newblock Chaotic weighted shifts on directed trees, 2024.
\newblock Preprint, arXiv:2303.03980.

\bibitem{GrossePeris}
Karl-G. Grosse-Erdmann and Alfredo Peris~Manguillot.
\newblock {\em Linear chaos}.
\newblock Universitext. Springer, London, 2011.

\bibitem{Hedlund}
J.~H. Hedlund.
\newblock Expansive automorphisms of {B}anach spaces. {II}.
\newblock {\em Pacific J. Math.}, 36(3):671--675, 1971.

\bibitem{Jablonski}
Zenon~Jan Jab{\l}o{\'n}ski, Il~Bong Jung, and Jan Stochel.
\newblock Weighted shifts on directed trees.
\newblock {\em Mem. Amer. Math. Soc.}, 216(1017):viii+106, 2012.

\bibitem{LopezPapachainrecontrees}
Antoni L\'opez-Mart\'inez and Dimitris Papathanasiou.
\newblock Shifts on trees versus classical shifts in chain recurrence.
\newblock {\em J. Differential Equations}, 433:Paper No. 113230, 25, 2025.

\bibitem{Martinez}
Rub\'en~A. Mart\'inez-Avenda\~no.
\newblock Hypercyclicity of shifts on weighted {$L^p$} spaces of directed trees.
\newblock {\em J. Math. Anal. Appl.}, 446(1):823--842, 2017.

\bibitem{MortensenWalmsley}
Andrew Mortensen and David Walmsley.
\newblock Chain recurrent shifts on trees, 2026.
\newblock Preprint, arXiv:2603.18878.

\end{thebibliography}

\end{document}